\documentstyle[leqno,11pt]{article}

\makeatletter
\@addtoreset{equation}{section}
\makeatother

\newcommand{\R}{I\hspace{-1.5mm}R}
\newcommand{\scriptR}{I\hspace{-0.9mm}R}
\newcommand{\N}{I\hspace{-1.4mm}N}

\newtheorem{Theorem}{Theorem}[section]

\newtheorem{Proposition}[Theorem]{Proposition}
\newtheorem{Lemma}[Theorem]{Lemma}
\newtheorem{Corollary}[Theorem]{Corollary}
\newtheorem{Remark}[Theorem]{Remark}
\date{}
\title{Value Functions for Bolza Problems with Discontinuous Lagrangians
and
Hamilton-Jacobi Inequalities\thanks{
Supported by the Human Capital and Mobility
Programme of
the European Union, contract number CHRX-CT94-0431}}
\author{Gianni Dal Maso\thanks{SISSA, via Beirut 2, 34014 Trieste, Italy
} \hspace{.1cm}
 and H\'el\`ene Frankowska\thanks{CNRS, ERS2064, Centre de
Recherche Viabilit\'{e}, Jeux, Contr\^{o}le, Universit\'{e} de
Paris-Dauphine, 75775 Paris Cx 16, France} \hspace{.1cm} }

\begin{document}
\maketitle
\noindent
{\bf Abstract.}

\noindent
We investigate the value function of the Bolza problem of the
Calculus of Variations
\begin{displaymath}
 V (t,x)=\inf \left\{ \int_{0}^{t} L (y (s),y' (s))ds \;+\;
\varphi (y(t)) : y \in W^{1,1} (0,t;\R^n),\; y(0)=x \right\},
\end{displaymath}
with a lower semicontinuous Lagrangian $L$ and a final cost $ \varphi $,
and
show that it is locally Lipschitz for $t>0$
whenever $L$ is locally bounded. It also satisfies
Hamilton-Jacobi inequalities in a  generalized sense.

When the Lagrangian is continuous, then the value function is the
unique lower semicontinuous solution
to the corresponding Hamilton-Jacobi equation, while for discontinuous
Lagrangian we characterize the value function by using the so
called contingent inequalities.

\vspace{ 5 mm}

\noindent
{\bf Key words.} Discontinuous Lagrangians,
 Hamilton-Jacobi equations, viability theory, viscosity solutions.

\vspace{ 5 mm}

\noindent
{\bf AMS-MOS subject classification:} 49L20 (primary),
49L25 (secondary).

\section{Introduction}
Consider the autonomous Bolza minimization problem of the Calculus of
Variations
\begin{equation}\label{B}
\inf \left\{ \int_{0}^{T} L (y (s),y' (s))ds \;+\;
\varphi (y(T)) ; y \in W^{1,1} (0,T;\R^n),\; y(0)=x \right\},
\end{equation}
where $ \varphi \colon  \R^n \mapsto \R_+ \cup  \{ + \infty \}$ and
$L\colon  \R^n
\times \R^n \mapsto \R_+$ are lower semicontinuous and $L (x, \cdot )$
is convex. The classical Lagrange problem (with the fixed final condition
$y (T)=x_T$) may be reduced to the above one by simply setting $ \varphi
(x_T)=0$ and $ \varphi = + \infty $ elsewhere.
Lagrange problems were studied in the literature from various points of
view.
If a Tonelli type
coercivity assumption holds true
\begin{equation} \label{T}
\exists \;  \Theta \colon  \R^n \mapsto \R_+ ,\;\;
 \lim_{|u| \rightarrow \infty } \frac{ \Theta (u)}{|u|} =+
\infty , \;\; \forall \; (x,u),\;\; L (x,u) \geq \Theta ( u ) ,
\end{equation}
then absolutely continuous minimizers do exist (see, for instance,
\cite[Cesari]{ces83}). Actually the very same condition yields also lower
semicontinuity of the value function associated to problem (\ref{B}).
It was observed in  \cite[Ambrosio, Ascenzi and Buttazzo]{ambr} that,
without
additional (boundedness) assumptions, minimizers are in general not
Lipschitz. This creates a major difficulty in developing the
Hamilton-Jacobi theory in a general case.
It was also shown in \cite{ambr} that minimizers  are Lipschitz
continuous,
under the local boundedness condition
\begin{displaymath}
 \forall \; R>0, \; \exists \; M>0,\; \exists \;  r>0,\;\; \forall \;
(x,u)
\in B (0,R) \times  B (0,r),\;\; L (x,u) \leq M ,
 \end{displaymath}
where $B(x,r)$ denotes the closed ball with centre $x$ and radius $r$.
The above condition may be rewritten
in the following equivalent way:
\begin{equation} \label{bound}
\ \ \ \ \ \ \
\forall \; x_0 \in \R^n, \;\exists \; M>0,\; \exists \;  r>0,\;\; \forall
\;
(x,u) \in B (x_0,r) \times  B (0,r),\;\; L (x,u) \leq M .
\end{equation}
Since any minimizer of problem  (\ref{B}) solves also an associated
Lagrange
problem, the same statement remains valid also for the Bolza problem,
whenever the infimum in (\ref{B}) is finite.

Consider the Hamilton-Jacobi equation
\begin{equation}\label{H-J}
\cases{V_t+H(x,-V_x)=0& in $\R_+^{\star} \times \R^n$,
\cr\cr
V (0, \cdot )= \varphi& in $\R^n$,
\cr}
\end{equation}
where $\R_+^{ \star }=\R_+ \backslash \{0\}$ and $H$ is the Hamiltonian
associated with $L$, defined by
\begin{equation}\label{H}
 H (x,p)= \sup_{u \in \scriptR^n} \left( \left\langle   p,u\right\rangle
 -L(x,u) \right),
\end{equation}
i.e., $H (x, \cdot )$ is the Legendre-Fenchel transform of $L (x, \cdot
)$.

It is a well known fact that  (\ref{H-J}) does not have smooth solutions
even when  $H$ and $\varphi $ are smooth, and one has to use generalized
solutions to treat first order equations of the above type. To prove
the uniqueness of the solution to  (\ref{H-J}) in the viscosity sense (as
done in \cite[Crandall and Lions]{cl2}), the
authors assumed, among other hypotheses, the continuity of $H$ and looked
for
bounded uniformly continuous solutions (in the viscosity sense).

The link between (\ref{H-J}) and (\ref{B}) is the value function
$V \colon  \R_+ \times  \R^n \mapsto \R_+ \cup \{+ \infty \}$,
defined by
\begin{displaymath}
 V (t,x)= \inf \left\{ \int_{0}^{t} L (y (s),y' (s))ds \;+\; \varphi
(y (t)) : y\in W^{1,1} (0,t;\R^n),\; y(0)=x \right\},
\end{displaymath}
which is, under some regularity assumptions, a viscosity solution to
(\ref{H-J}), i.e., it is a
supersolution (resp.  subsolution) with derivatives replaced by
subdifferentials (resp. superdifferentials). The above definition of value
function is somewhat different from the usual one,
but, our problem being autonomous,
we have found it more convenient for our purposes.
The reader accustomed to
different definitions can easily make suitable changes to derive
similar results for other value functions.

There exists an interplay between subdifferentials of the value function
and
minimizers. For instance, in Section 3 we show that, if
$ y$ is a minimizer to the last problem and if  the subdifferential of
$V$ at $ (t,x)$ is nonempty, then the difference quotients
$ \left\{  \frac{ y (h)-x}{h} \right\}$ are bounded
(Proposition~\ref{cont}).

Many results about  Bolza problems with smooth data were extended to the
case of lower
semicontinuous  $L$ and $\varphi $. In such case, under Tonelli's
assumption
 (\ref{T}), $H $
is merely upper semicontinuous and, for this reason, it is natural to look
for an extension of the viscosity solutions theory to this situation. Also
in such general case $V$ is only lower semicontinuous, which creates
additional difficulties in formulating maximum principles yielding
uniqueness. However we prove that $V$ is locally Lipschitz on
$\R_+^{\star } \times \R^n$,
whenever $L$ is locally bounded, even if the data $L$ and
$\varphi $ are discontinuous (Corollary~\ref{vlip}).

In the case of lower semicontinuous solutions and
Hamilton-Jacobi-Bellman equations related to Mayer's problem of optimal
control theory, a maximum principle for lower semicontinuous functions was
proposed in \cite[Barron and Jensen]{Barron-Jensen1} to deduce uniqueness
of
solutions to the  Hamilton-Jacobi equation (corresponding to Mayer's
problem).
In \cite{f91hjb,f93hjb} the same uniqueness result was obtained
by exploiting properties of the epigraph of the solutions of the
Hamilton-Jacobi equation.

When $L$ is continuous, we prove (Theorem~\ref{smaller}) that $V$ is
the unique lower semicontinuous function which satisfies the initial
condition $V(0, \cdot )= \varphi$ and solves the Hamilton-Jacobi
equation (\ref{H-J}) in the following
sense:
\begin{equation} \label{cd220}
 \forall \;(t,x) \in \mbox{\rm dom}(V), \;t >0,\; \exists \; u
 \in\R^{n},\; \;
 D_{\uparrow }V (t,x) (-1, u) \leq -L (x,  u) ,
\end{equation}
 \begin{equation} \label{cd210}
 \forall \; (t,x) \in \mbox{\rm dom}(V), \;
 \forall \; u \in \R^n,\;\; D _{\downarrow }V (t,x) (1,-u) \leq L (x,u) .
\end{equation}
where $D _{\uparrow }V (t,x)$ and $D _{\downarrow }V (t,x)$ are the
lower and upper contingent
derivatives of $V$ at $ (t,x)$, whose definition is recalled in Section 2.
We underline that they coincide with Dini's lower and upper derivatives
where $V$ is locally Lipschitz.

Contingent inequalities for discontinuous functions
were introduced in \cite[Aubin]{aublyap}
to study lower semicontinuous Lyapunov functions. They were introduced
independently in the context of Lipschitz functions (and Dini's
directional
derivatives)
in  \cite[Subbotin]{subbo}
to investigate Isaacs' equation of differential games. In
\cite{f87hjc} contingent inequalities were used to study lower
semicontinuous supersolutions of
Hamilton-Jacobi-Bellman equation of optimal
control.

Under the same assumptions we prove also (Theorem~\ref{hju1}) that
the restriction of $V$ to $\R_+^{ \star } \times \R^n$ is the
unique locally Lipschitz function which satisfies the initial
condition
\begin{equation}
 \label{boundary}
 \liminf_{h \rightarrow 0+, \;y \rightarrow x} V (h,y) = \varphi  (x)
\end{equation}
and solves the Hamilton-Jacobi
equation (\ref{H-J})
in the following generalized sense (which is weaker than
(\ref{cd220}) and (\ref{cd210})):
\begin{displaymath}
\forall \;(t,x) \in \R_+^{ \star }\times \R^n,\;\; \forall \;
(p_t,p_x)
\in \partial _- V (t,x),\;\; p_t+ H (x,-p_x)=0 ,
\end{displaymath}
where $ \partial _-V (t,x)$ denotes the subdifferential of $V$ at $
(t,x)$.

Moreover we prove that the restriction of $V$ to $\R_+^{ \star }
\times \R^n$ is the
unique locally Lipschitz function which satisfies the initial
condition (\ref{boundary}) and solves the Hamilton-Jacobi
equation (\ref{H-J})
in the viscosity sense:
\begin{displaymath}
 \forall \; (t,x) \in \R_+^{ \star }\times \R^n, \;\; \forall \;
(p_t,p_x)
\in \partial _- V (t,x),\;\; p_t+ H (x,-p_x) \geq 0 ,
\end{displaymath}
\begin{displaymath}
  \forall \; (t,x) \in \R_+^{ \star }\times \R^n, \;\; \forall \;
(p_t,p_x)
\in \partial _+ V (t,x),\;\; p_t+ H (x,-p_x) \leq 0 ,
\end{displaymath}
where $ \partial _+V (t,x)$ denotes the superdifferential of $V$ at $
(t,x)$. We have been not able to prove that there is only one
lower semicontinuous viscosity solution. For this reason we
have to adopt a rather unusual notion of solution in our uniqueness
result for lower semicontinuous solutions.

When $L$ is discontinuous and locally bounded,
we prove (Theorem~\ref {last}) that $V$ is the
unique lower semicontinuous  function
which satisfies the initial condition $V(0, \cdot )= \varphi$
and solves the Hamilton-Jacobi
equation (\ref{H-J}) in the following
sense:
\begin{displaymath}
 \forall \; (t,x) \in \mbox{\rm dom}(V),\; t>0, \;
 \exists\; u\in\R^{n}, \;\;
 D _{\uparrow }V (t,x) (-1,u) \leq -L (x,u) ,
\end{displaymath}
\begin{displaymath}
\forall \; (t,x) \in \mbox{\rm dom}(V),\;
\forall \; u\in\R^{n},\;\; D_{\downarrow } V (t,x) (1,-u) \leq
L^+(x,u) ,
\end{displaymath}
where $L^{+}$ is  defined by
\begin{displaymath}
L^+(x,u) = \limsup_{h \rightarrow 0+} \frac{1}{h} \inf \left\{
\int_{-h}^{0} L (y (s),y' (s))ds  : y (-h)=x-hu ,\; y (0)=x\right\}
\end{displaymath}
(see \cite{dalm} and \cite[Amar, Bellettini, Venturini]{amar}).

Under the same hypotheses we prove also (Theorem~\ref{ddd050}) that
the restriction of $V$ to $\R_+^{ \star } \times \R^n$ is
the
unique locally Lipschitz function which satisfies the initial
condition (\ref{boundary})
together with the additional condition
\begin{equation} \label{ddd15}
\forall \; x\in \R^{n}, \; \; \forall \; \lambda>0, \;\;
\lim_{\scriptstyle h \rightarrow 0+, \;y \rightarrow x
\atop \scriptstyle |y-x| \le \lambda h} V (h,y) = \varphi  (x),
\end{equation}
and solves the Hamilton-Jacobi
equation (\ref{H-J}) in the
following sense:
\begin{displaymath}
\forall \;(t,x) \in \R_+^{ \star }\times \R^n,\;\; \forall \;
(p_t,p_x)
\in \partial _- V (t,x),\;\; p_t+ H (x,-p_x) \geq 0,
\end{displaymath}
\begin{displaymath}
\forall \;(t,x) \in \R_+^{ \star }\times \R^n,\;\; \forall \;u \in
\R^n,
\;\; D_{\downarrow }V (t,x) (1,-u) \leq L^+ (x,u).
\end{displaymath}
In all these theorems, the uniqueness is obtained by proving
suitable comparison results for the corresponding notions of sub- and
supersolution. In particular, we extend here a result from \cite{f89hjc}
proved in the context
of Mayer's problem with bounded dynamics,
and show (Corollary~\ref{gg1}) that, if (\ref{bound}) holds true, then the
value function is the smallest lower semicontinuous function satisfying
the initial condition $ V (0, \cdot )= \varphi $ and the contingent inequality
(\ref{cd220}).

Recently solutions to the Hamilton-Jacobi equation of a nonautonomous
Bolza
problem  were investigated in \cite[Galbraith]{galbraith}. However the
results of
\cite{galbraith} do not overlap with ours, since the assumptions
of that article imply that the Hamiltonian is locally Lipschitz.

\section{Preliminaries}

Let $K\subset \R^n$ be a nonempty subset and $x\in K$. The
contingent cone $T_K(x)$ to $K$ at $x$ is defined by
\[
v\in T_K(x) \Longleftrightarrow \liminf_{h\to 0+}
\frac{\mbox{\rm dist}(x+hv,K)}{h}=0 .
\]
The negative polar cone $T^-$ to a subset $T\subset \R^n$ is given by
\[
T^-=\{v\in \R^n\,:\;\forall \; w\in T\,,\;\langle v,w \rangle\leq 0\} .
\]
We recall the following result due to Cornet \cite{Cornet} (see also
\cite{af90sva}).
\begin{Theorem} \label{corn}
Let $K\subset \R^n$ be a closed subset and $x\in K$. Then
\begin{displaymath}
\mathop{\mbox{
{\rm Liminf}}}_{ \textstyle{y \rightarrow x \atop y \in K}}
\overline{co} \; T_K (y) = C_K (x) \subset T_K (x),
\end{displaymath}
where  $\mbox{ {\rm Liminf}}$ denotes the topological lower limit (in the
Painlev\'e-Kuratowski sense) and
$C_K (x)$ denotes Clarke's tangent cone to $K$ at $x$.
\end{Theorem}

Let $W\colon  \R^n\to \R \cup \{+ \infty \}$ be a lower semicontinuous
function.
The subdifferential of $W$ at $x\in \mbox{\rm dom}(W)$  is defined by
\[
\partial_- W(x)=\left\{
p\in \R^n\;:\; \liminf_{y\to x}\frac{ W(y)-
W(x)-\langle p,y-x\rangle}{|y-x|} \geq 0\right\} .
\]
By  \cite{f89hjc} (or \cite{af90sva})
\begin{equation} \label{subd}
 p\in\partial_-
W(x) \;\;\Longleftrightarrow\;\; (p,-1)\in\left[T_{{\cal E}pi(
W)}(x,W(x))\right]^- ,
\end{equation}
where ${\cal E}pi(W)$ denotes the epigraph of $W$ defined by
\begin{displaymath}
 {\cal E}pi(W) := \left\{ (x,r)\in \R^n\times\R : r \geq W (x) \right\} .
\end{displaymath}

An equivalent definition of subdifferental uses directional derivatives of
$W$ defined by
\begin{equation}\label{dirder}
\forall \; u \in \R^n,\;\; D _{\uparrow } W (x) (u) =
\liminf_{ \textstyle {h \rightarrow 0+ \atop v \rightarrow
u}} \frac{W (x+hv) - W (x)}{h} \;\;.
\end{equation}
Clearly for every $x\in \mbox{\rm dom}(W)$
\begin{equation}
\label{star}
 {\cal E}pi (D _{\uparrow } W (x) ( \cdot ) ) =
 T_{{\cal E}pi (W)} (x,W (x)) ,
\end{equation}
and therefore
\begin{equation}
 \label{subd1}
 p \in \partial _- W (x)  \;\; \Longleftrightarrow \;\; \forall \; v\in
\R^n,
\;\; \left\langle   p,v\right\rangle \leq D _{\uparrow } W (x) (v) .
\end{equation}

The upper directional derivative of $W$ at $x$ in the direction $u$ is
 defined by
\begin{displaymath}
\forall \; u \in \R^n,\;\; D _{\downarrow }  W (x) (u) =
\limsup_{ \textstyle {h \rightarrow 0+ \atop v \rightarrow
u}} \frac{W (x+hv) - W (x)}{h} .
\end{displaymath}
The superdifferential $ \partial_+ W (x)$ of $W$ at $x$ is defined by
$ \partial_+ W (x)=- \partial_- (-W) (x)$ or, equivalently, by
\begin{equation}
 \label{supd1}
 p \in \partial_+ W (x)  \;\; \Longleftrightarrow \;\; \forall \; v\in
\R^n,
\;\; \left\langle   p,v\right\rangle \geq D_{\downarrow } W (x) (v) .
\end{equation}

We shall
need the following
version of Rockafellar's result (see \cite{Rockafellar}).
\begin{Lemma}\label{Rockafellar}
Let $x\in \mbox{\rm dom}(W)$ and let
$(p,0)\in\left[T_{{\cal E}pi(W)}(x,W(x))\right]^-$ be such that
$p\neq 0$. Then there exist
$x_\varepsilon $ converging to $x$ (as $ \varepsilon  \rightarrow 0+$) and
\begin{displaymath}
 (p_\varepsilon,q_\varepsilon)\in\left[T
_{{\cal E}pi(W)}(x_\varepsilon,W(x_\varepsilon ))\right]^-
\end{displaymath}
converging to $ (p,0)$ as  $ \varepsilon  \rightarrow 0+$
such that $q_{ \varepsilon }<0$.
\end{Lemma}

A closed subset $K$ of $\R^n$ is called a viability domain
of a set-valued map $G\colon \R^n \leadsto  \R^n$ if for every $x\in K$
\[
G(x)\cap T_K(x)\neq\emptyset \;\;.
\]
The following formulation summarizes several versions of the viability
theorem
(see \cite{Aubin1} and  \cite{Aubin}).
\begin{Theorem}[Viability]\label{viability}
Suppose that $G\colon \R^n\leadsto \R^n$ is an upper semicontinuous
set-valued map with compact convex values. For a
closed subset
$K\subset \R^n$ the following conditions are equivalent:
\begin{description}
\item[\mbox{\rm (a)}] $K$ is a viability domain of $G$;
\item[\mbox{\rm (b)}]  $G(x)\cap\overline {co}\,T_K(x)\neq\emptyset$
for every $x\in K$;
\item[\mbox{\rm (c)}]  for every $x_0\in K$ there exist
$ \varepsilon >0$ and a solution $x\colon {[0,\varepsilon[}\mapsto K$
to the Cauchy problem
\end{description}
\begin{equation}\label{diG}
\cases{
x'(t)\in G(x(t)) ,
\cr
x(0)=x_0 .&
\cr}
\end{equation}
\end{Theorem}

The equivalence (a) $\Longleftrightarrow$ (b)
was proved in \cite{Ushakov}. This proof was simplified in
\cite[p.85]{Aubin}. The fact that (a)
$\Longleftrightarrow$ (c) was first proved by
Bebernes and Schuur in \cite{beb}. A proof can be found
in \cite{AubinCellina} or \cite{Aubin1}, \cite{Aubin}.

\section{The value function of the Bolza problem}

Consider a lower semicontinuous function
$L\colon  \R^n \times \R^n \mapsto \R_+ $
and an extended lower semicontinuous function $ \varphi \colon  \R^n
\mapsto \R_+ \cup \left\{ + \infty  \right\}$, not identically equal to
$+\infty$.
        Throughout the whole paper we suppose  that $L$ is convex in the
second variable and that the
coercivity assumption  (\ref{T}) holds true.
Without any loss of the generality we also  assume that  $ \Theta$ is
convex.

\vspace{ 5 mm}

Given $t_0>0$ and $x_0 \in \R^n$, let us consider the Bolza problem:
\begin{displaymath}
\mbox{ {\rm minimize}} \int_{0}^{t_0} L (y(s), y' (s))ds +
\varphi  (y(t_0))
\end{displaymath}
over all absolutely continuous functions $y \in W^{1,1} (0,t_0; \R^n)$
satisfying
the initial condition $y (0)= x_0$.
The dynamic programming approach associates with this problem the family
of
problems ($t\geq 0$, $x\in\R^n$):
\begin{displaymath}
\mbox{ {\rm minimize}} \int_{0}^{t} L (y (s), y' (s))ds +
\varphi  (y(t))
\end{displaymath}
over all absolutely continuous functions $y \in W^{1,1} (0,t; \R^n)$
satisfying
$y (0)= x$. The value function $V\colon  \R_+ \times \R^n \mapsto  \R_+
\cup \left\{ + \infty  \right\}$ is defined by
\begin{equation} \label{value}
 V (t,x)= \inf \left\{\int_{0}^{t} L (y (s), y' (s))ds + \varphi  (y
(t)) :  y (0)=x \right\} .
\end{equation}

 \begin{Proposition} \label{lsc} Under the above assumptions
for all $ (t,x) \in \R_+ \times \R^n$ the infimum in
(\ref{value}) is attained (it may be infinite)
and $V$ is lower semicontinuous on $\R_+ \times \R^n$.
Furthermore, if $L$ is locally bounded, then $V$ has finite
values on $\R_+^{\star} \times \R^n$ and
satisfies (\ref{boundary}) and (\ref{ddd15}) .
\end{Proposition}

{\bf Proof} --- \hspace{ 2 mm}
The existence of a minimizer is a well known result. For the reader's
convenience, we sketch the proof of the
lower semicontinuity of $V$.
Consider a sequence $ (t_i,x_i) \in \R_+ \times \R^n$ converging to $
(t,x)$ such that $V (t_i,x_i)$ converge to $\liminf_{ (s,y) \rightarrow
(t,x)}V (s,y)$.  Let $ y_i$ be the corresponding minimizers of
(\ref{value}) with $ (t,x)$ replaced by $ (t_i,x_i)$. If $\displaystyle
\lim_{i \rightarrow + \infty }V (t_i,x_i) = + \infty $, then
$$V (t,x)
\leq \lim_{i \rightarrow + \infty }V (t_i,x_i)\;. $$
Assume next that the above limit is finite.
Hence for some $M$ for all  $i$,
\begin{displaymath}
 \int_{0}^{t_i}L ( y_i (s), y_i ' (s))ds \leq M .
\end{displaymath}
Set $y_i(s)=x_i$ for $s \geq t_i$.
By the coercivity assumption (\ref{T}), the sequence
$y'_i$ is
equiintegrable on $[0,t+1]$.
This, the Dunford-Pettis criterion (see, e.g., \cite[Theorem
II.25]{mayer}),
and the Ascoli-Arzel\`a
theorem yield the existence of a
subsequence $y_{i_k}$ converging uniformly to some $y$ on $[0,t+1]$ such
that
$y_{i_k}'$ converge to $ y'$ weakly in $L^1$. We denote
this subsequence again by $y_i$.
Then $y (0)=x$ and $y_i (t_i)\rightarrow y(t)$.
Fix $0< \varepsilon <t$. Then for all large $i$,
\begin{displaymath}
 \int_{0}^{t_i}L (y_i (s), y_i ' (s)))ds
\geq  \int_{0}^{t- \varepsilon }L (y_i (s), y_i ' (s)))ds .
\end{displaymath}
Since $L$ is lower
semicontinuous and convex in the second variable,
from the lower semicontinuity theorems by Olech \cite{Olech} and Ioffe
\cite{Ioffe} (see also \cite[Theorem 2.3.3]{butt}) it follows  that
\begin{displaymath}
\liminf_{i \rightarrow \infty }\int_{0}^{t- \varepsilon }L (y_i (s), y_i '
(s)))ds \geq \int_{0}^{t- \varepsilon }
L (y (s), y ' (s))ds .
\end{displaymath}
Thus
\begin{displaymath}
\liminf_{i \rightarrow \infty } \int_{0 }^{t_i}L (y_i (s), y_i ' (s)))ds
\geq \int_{0}^ {t- \varepsilon }L (y (s), y ' (s))ds .
\end{displaymath}
Since $\int_{0}^{t- \varepsilon }L (y (s), y ' (s))ds$ converges to
$ \int_{0}^{t}L (y (s), y ' (s))ds$ when $\varepsilon
\rightarrow 0+$, the
 lower semicontinuity of  $ \varphi $ yields lower semicontinuity of  $V$.

Since $ \varphi $ is not identically $+ \infty $, it
is clear that $V$ is finite on $\R_+^{\star} \times \R^n$
whenever $L$ is locally bounded.

To prove (\ref{ddd15})
consider a sequence $ (h_i,x_i) \rightarrow (0+,x)$ such that $
(x_i-x)/h_i$ is bounded. Then, by lower semicontinuity of $V$,
\begin{displaymath}
 \liminf_{i \rightarrow  \infty }V (h_i,x_i) \geq V (0,x)= \varphi  (x),
\end{displaymath}
\begin{displaymath}
 V (h_i,x_i) \leq  \varphi (x) + \int_{0}^{h_i}L (x_i +s (x-x_i)/h_i,
(x-x_i)/h_i)ds .
\end{displaymath}
Passing to the upper limit when $i \rightarrow  \infty $, and using
the fact that
$L$ is locally bounded, we obtain (\ref{ddd15}).
Condition (\ref{boundary}) follows now from (\ref{ddd15}) and
from the lower semicontinuity of $V$.
$\; \; \Box$
\begin{Proposition} \label{cont}
 Let $ (t,x) \in \mbox{\rm dom}(V)$, with $t>0$, and let $y$ be a
minimizer of  (\ref{value}).
Assume either that
$ \partial _- V (t,x)$ is nonempty or that (\ref{bound}) holds.
Then the set
\begin{equation} \label{dm001}
 \left\{ \frac{y (h) -x }{h} \right\}_{h>0}
\end{equation}
is bounded in $\R^n$.
Furthermore, if
\begin{equation} \label{dm15'}
\frac{y (h_i)-x}{h_i} \rightarrow u
\end{equation}
for some $h_i \rightarrow 0+$, then
\begin{equation} \label{coin}
 D _{\uparrow }V (t,x) (-1,u) \leq -L (x,u) .
\end{equation}
Consequently, $V$ is a viscosity supersolution to (\ref{H-J}), i.e.,
\begin{equation} \label{supersol}
 \forall \; (t,x) \in \R_+^{ \star } \times \R^n ,\;\; \forall \;
(p_t,p_x)
\in
\partial _-V (t,x),\;\;p_t+ H (x,-p_x) \geq 0 .
\end{equation}
\end{Proposition}

\vspace{ 3 mm}

{\bf Proof} --- \hspace{ 2 mm}
First we observe that for all $ 0 \leq h \leq t$
\begin{equation}
 \label{equa}
 V (t,x)= V (t-h, y (h))+ \int_{0}^{h}L (y (s), y' (s))ds .
\end{equation}

{\bf Step 1.}
Let us prove that (\ref{dm001}) is bounded.
If (\ref{bound}) holds, then (\ref{dm001}) follows from \cite{ambr}.
Suppose now that $ \partial _- V (t,x)$ is nonempty. If
(\ref{dm001}) does not hold,
then there
exists a sequence $h_i \rightarrow 0+$ such that
\begin{equation}
 \label{dm002}
\frac{|  y (h_i)-  x |}{h_i} \rightarrow  +\infty .
\end{equation}
Taking a subsequence, still denoted by $h_i$,
we may assume that for some $v \in \R^n$
\begin{displaymath}
 \frac{y (h_i)- x}{| y (h_i)- x |}
\rightarrow v .
\end{displaymath}
Fix  $ (p_t,p_x) \in \partial _- V (t,x)$.
By the definition of subdifferential
\begin{equation}
 \label{dm18''}
 \liminf_{h \rightarrow 0+} \frac{V (t-h, y (h))- V (t,x) -
\left\langle   (p_t,p_x), (-h, y (h)-x)\right\rangle }{h + |
y (h)- x |} \geq 0 ,
\end{equation}
which yields, by (\ref{equa}),
\begin{displaymath}
 \limsup_{i \rightarrow \infty }
 \frac{1}{h_i + | y (h_i)- x |} \int_{0}^{h_i}L (y (s), y'(s))ds
\leq -\left\langle  p_x,v\right\rangle .
\end{displaymath}
By (\ref{T}) this inequality implies that there exists a constant
$c>0$ such that
\begin{displaymath}
\forall \; i\in\N, \;\;
 \frac{1}{h_i} \int_{0}^{h_i} \Theta(y'(s))ds \leq
 c\,\Big( 1 + \frac{|  y (h_i)-  x |}{h_i} \Big).
\end{displaymath}
By Jensen's inequality we obtain
\begin{displaymath}
\forall \; i\in\N, \;\;
\Theta\Big(\frac{  y (h_i)-  x }{h_i}\Big) \leq
 c\,\Big( 1 + \frac{|  y (h_i)-  x |}{h_i} \Big) .
\end{displaymath}
By (\ref {dm002}) this contradicts
the coercivity assumption (\ref{T}) and ends the proof of
our claim.

\vspace{ 5 mm}

{\bf Step 2.} Let us fix $R>0$ and $\varepsilon >0$.
We want to prove that there exists $\delta>0$ such that
\begin{equation}
\label{dm18'}
\forall \; \xi\in B(x,\delta), \;\; \forall \; u\in B(0,R), \;\;
L(\xi,u) \geq L(x,u) - \varepsilon .
\end{equation}
We start by observing that the function $L(x,\cdot)$ is continuous,
since it is convex and finite valued.
Thus for every $v\in{B(0,R)}$ there exists
$\rho_1=\rho_1(v,\varepsilon)>0$ such that
\begin{displaymath}
\forall\; u\in B(v,\rho_1),\;\; L(x,v)\ge
L(x,u) -  \frac{\varepsilon}{2} .
\end{displaymath}
As $L$ is lower semicontinuous, for every $v\in{B(0,R)}$
there exists
$\rho_2=\rho_2(v,\varepsilon)>0$ such that
\begin{displaymath}
\forall\; \xi\in B(x,\rho_2),\;\;  \forall\; u\in B(v,\rho_2),\;\;
L(\xi,u)\ge  L(x,v) -  \frac{\varepsilon}{2} .
\end{displaymath}
Putting $r=r(v,\varepsilon)=
\min\{\rho_1(v,\varepsilon), \rho_2(v,\varepsilon)\}$, it follows that
\begin{displaymath}
\forall\; \xi\in B(x,r),\;\;  \forall\; u\in B(v,r),\;\;
L(\xi,u)\ge  L(x,u) -  \varepsilon .
\end{displaymath}
By compactness there exists a finite set $\{v_1,\ldots,v_k\}$ contained in
${B(0,R)}$ such that ${B(0,R)}\subset B(v_1,r_1) \cup
\cdots \cup B(v_k,r_k)$, where $r_i=r(v_i,\varepsilon)$.
It is then clear that (\ref{dm18'}) is satisfied
with $\delta=\min\{r_1,\ldots,r_k\}$.

\vspace{ 5 mm}

{\bf Step 3.}
Consider now a sequence $h_i \rightarrow 0+$ and $u\in\R^n$ such that
(\ref{dm15'}) holds. By (\ref{equa}) we have
\begin{equation} \label{dm19}
D _{\uparrow }V (t,x) (-1,u) \leq -\limsup_{i \rightarrow  \infty }
\frac{1}{h_i} \int_{0}^{h_i} L (y (s), y' (s))ds .
\end{equation}
{}From (\ref{equa}) and (\ref{dm18''}) we obtain
\begin{displaymath}
\limsup_{i \rightarrow  \infty }
\frac{1}{h_i} \int_{0}^{h_i} L (y (s), y' (s))ds
\leq p_t - \langle p_x,u \rangle .
\end{displaymath}
 By (\ref{T}) this implies that there exists a constant $C>0$ such that
\begin{equation} \label{dm00}
\forall \; i\in\N, \;\; \frac{1}{h_i} \int_{0}^{h_i} \Theta ( y'
(s))ds\leq C .
\end{equation}

Let us fix $R>0$ and $\varepsilon >0$. By (\ref{dm18'}) for $i$ large
enough we have
\begin{equation} \label{dm01}
L (y (s), y' (s)) \geq L (x, y' (s)) - \varepsilon
\end{equation}
for every $s\in[0,h_i]$ such that $|y' (s)|<R$. For every $i$ let
\begin{displaymath}
A_i=\{ s\in[0,h_i] : |y' (s)| < R\},\qquad
B_i=\{ s\in[0,h_i] : |y' (s)|\geq R\},
\end{displaymath}
and let $\theta(R)=\min\{\Theta(v) : |v|\geq R \}$. Then
\begin{displaymath}
|B_i|\leq \frac{1}{\theta(R)} \int_{0}^{h_i} \Theta ( y' (s))ds
\leq \frac{C}{\theta(R)} h_i ,
\end{displaymath}
where $|\cdot |$ denotes the Lebesgue measure. Consequently
$|A_i|\geq \tau(R)\,h_i$, where
\begin{displaymath}
\tau(R) = \Big(1 - \frac{C}{\theta(R)} \Big) \rightarrow 1
\ \mbox{  as  } \ R \rightarrow +\infty .
\end{displaymath}
By (\ref{dm01}) we have
\begin{displaymath}
\frac{1}{h_i} \int_{0}^{h_i} L (y (s), y' (s))ds \geq
\frac{\tau(R)}{|A_i|}
\int_{A_i} L (x, y' (s))ds - \varepsilon,
\end{displaymath}
and by Jensen's inequality we obtain
\begin{equation} \label{dm02}
\frac{1}{h_i} \int_{0}^{h_i} L (y (s), y' (s))ds \geq
\tau(R) \, L (x, u_i) - \varepsilon ,
\end{equation}
where $u_i=|A_i|^{-1} \int_{A_i} y' (s)ds$.
For every $v\in \R^{n}$ let $\omega(v)=| L(x,u+v) - L(x,u) |$.
Since $L (x, \cdot )$
is continuous at $u$, the function $\omega$ is continuous
at $0$ and $\omega(0)=0$.
{}From (\ref{dm02}) we have
\begin{equation} \label{dm03}
\frac{1}{h_i} \int_{0}^{h_i} L (y (s), y' (s))ds \geq
\tau(R)  \, L (x, u)
- \omega(u_i-u) - \varepsilon .
\end{equation}

In order to estimate $|u_i-u|$, we notice that
\begin{displaymath}
|u_i-u| \leq \Bigl| \frac{1}{|A_i|} \int_{A_i} y' (s)ds
- \frac{1}{h_i} \int_{0}^{h_i} y' (s) ds \Bigr| + \varepsilon_i ,
\end{displaymath}
where
\begin{displaymath}
\varepsilon_i = \Bigl| \frac{y(h_i)-x}{h_i} - u \Bigr| \rightarrow 0
\ \mbox{  as  } \  i \rightarrow \infty .
\end{displaymath}
Therefore
\begin{equation} \label{dm008}
|u_i-u|\leq \frac{1}{|A_i|} \int_{B_i} |y' (s)| ds
+ \Big(\frac{1}{|A_i|} - \frac{1}{h_i} \Big)
 \int_{0}^{h_i} |y' (s)| ds + \varepsilon_i .
 \end{equation}
 Let us define
\begin{displaymath}
\zeta(R) = \min_{|v|\geq R} \frac{\Theta(v)}{|v|} .
\end{displaymath}
By (\ref{T}) $\zeta(R) \rightarrow +\infty$ as $R\rightarrow +\infty$,
and there exists a constant $\alpha\geq 0$ such that
$\Theta(v)\geq |v| - \alpha$ for every $v\in\R^n$. {}From (\ref{dm008})
we obtain
 \begin{displaymath}
 |u_i-u| \leq \Big(\frac{1}{\tau(R)\zeta(R)} + \frac{1}{\tau(R)} - 1\Big)
 \frac{1}{h_i} \int_{0}^{h_i} \Theta ( y' (s))ds + \alpha\,
 \Big(\frac{1}{\tau(R)} - 1\Big) + \varepsilon_i
 \end{displaymath}
 \begin{displaymath}
 \leq \frac{C}{\tau(R)\zeta(R)} + (C+\alpha)
 \Big(\frac{1}{\tau(R)} - 1\Big) + \varepsilon_i .
 \end{displaymath}
Since $\tau(R) \rightarrow 1$ and $\zeta(R) \rightarrow +\infty$
as $R\rightarrow +\infty$,  we have
\begin{displaymath}
|u_i-u| \leq \sigma(R) + \varepsilon_i ,
\end{displaymath}
where $\sigma(R)\rightarrow 0$ as $R\rightarrow +\infty$.
Therefore we deduce from (\ref{dm19}) and (\ref{dm03}) that
\begin{displaymath}
D _{\uparrow }V (t,x) (-1,u) \leq - \tau(R)
L(x,u) + \omega(\sigma(R)) + \varepsilon ,
\end{displaymath}
for every $R>0$ and $\varepsilon >0$. Taking the limit as
$R\rightarrow\infty$
and $\varepsilon\rightarrow 0+$ we obtain (\ref{coin}).

\vspace{ 5 mm}

{\bf Step 4.} Let us prove (\ref{supersol}).
Pick $ (p_t,p_x) \in \partial _-V (t,x)$ and $u$ satisfying
(\ref{dm15'}).
{}From (\ref{subd1}) and (\ref{coin}) we get
 $ p_t + \left\langle   -p_x,u\right\rangle -L (x,u) \geq 0$. The
 conclusion follows from the definition (\ref{H}) of $H$.
$\; \; \Box$

\begin{Proposition} \label{minlip}
Assume that $L$ is locally bounded and let $ (t_0,x_0) \in \R_+^{ \star }
\times \R^n$. Then  there
exist $r>0$ and $\delta >0$ such that for all $ (t,x) \in B ( (t_0,x_0),
\delta )$  every minimizer $ y ( \cdot; t, x)$ of  (\ref{value})  is
$r-$Lipschitz.
\end{Proposition}

{\bf Proof} --- \hspace{ 2 mm}
According to  \cite[Theorem 4.1]{ambr}
for every minimizer $y ( \cdot ):= y ( \cdot ;t,x)$ of
(\ref{value}) there exists $k (t,x)$
such that for some $p (s) \in \partial _{u} L (y (s),y' (s))$
\begin{displaymath}
 \left\langle   p (s),y' (s)\right\rangle -L (y (s),y' (s)) =k (t,x)\;\;
\mbox{ {\rm for a.e.\ }} s \in [ 0,t] .
\end{displaymath}
On the other hand, since $L$ is locally bounded, for some $M_{1}>0$ and for
all
$(t,x)$ sufficiently close to
$(t_0,x_0)$
\begin{displaymath}
 \int_{0}^{t} L (y (s;t,x), y' (s;t,x))ds \leq M_{1}.
\end{displaymath}
This and coercivity assumption (\ref{T}) imply that $y'(\cdot;t,x)$ are
equiintegrable and therefore $y(\cdot;t,x)$ is uniformly bounded in
$L^{\infty}$ for
$(t,x)$ near $(t_0,x_0)$.
Furthermore, there exists $R>0$ such that for all $ (t,x)$ sufficiently close to
$ (t_0,x_0)$
the sets
\begin{displaymath}
 A(t,x) := \left\{ s\in [0,t] : |y' (s;t,x)| \leq R \right\}
\end{displaymath} have positive measure.
Since $L$ is locally bounded,
$$\sup_{s \in [0,t]}
\sup_{u \in B (0,2R)}L (y(s;t,x), u)\le M_2 <
\infty $$
for all $ (t,x)$ sufficiently close to
$ (t_0,x_0)$.
This implies that for some $M_3<+\infty$
 $L(y (s;t,x), \cdot )$ is $M_3-$Lipschitz on
$B(0,R)$ for
$(t,x)$ near $(t_0,x_0)$ (see, e.g., \cite[Proof of Theorem 2.1]{aubconvex}).
Hence for almost every $s \in A(t,x)$ we have
$|\langle p(s), y'(s)\rangle| \le M_{3}R$, which implies
\begin{displaymath}
 |k(t,x)| \leq M_2 + M_3 R .
\end{displaymath}
Consequently,  $k(t,x)$ is bounded in a neighborhood of  $ (t_0,x_0)$.
By \cite[Proof of Theorem 4.2]{ambr}, $y ( \cdot ;t,x)$ are Lipschitz with
the same Lipschitz constant whenever $ (t,x)$ is sufficiently close to $
(t_0,x_0)$. $\; \; \Box$

\begin{Corollary} \label{vlip}
If $L$ is locally bounded, then
$V $ is locally Lipschitz on $\R_+^{ \star } \times \R^n$.
\end{Corollary}

{\bf Proof} --- \hspace{ 2 mm}
        Fix $ (t_0,x_0) \in \R_+^{ \star } \times \R^n \subset \mbox{\rm
        dom}(V)$.
By Proposition
\ref{minlip}, there exist $r>0$ and $\delta >0$ such that for all $ (t,x)
\in
B ( (t_0,x_0), \delta )$  every minimizer $ y ( \cdot; t, x)$ of
(\ref{value}) is $r-$Lipschitz. We may assume that $5\delta<t_0$.
Let $(t_1, x_1)$ and $(t_2, x_2)$ be two distinct points of
$B ( (t_0,x_0), \delta )$, let $h_1=|t_1-t_2| + |x_1-x_2|$, and
$s_1= h_1-t_1+t_2$. Let $u_1\in\R^n$ be
such that $y(s_1;t_2,x_2)= x_1+h_1u_1$. Then $0< h_1<t_1$,
$0\le s_1 \le 2h_1$, and
\begin{equation}\label{dm301}
|u_1| \leq \frac{|y(s_1;t_2,x_2)-x_2|}{h_1} +
\frac{|x_2-x_1|}{h_1} \leq 2r + 1.
\end{equation}
Let $y_1\colon [0,t_1]\mapsto\R^n$ be the function defined by
\begin{displaymath}
 y_1 (s)=\cases{x_1+su_1&if $0\le s\le h_1$,
 \cr\cr
 y(s-t_1+t_2; t_2,x_2)&if $h_1\le s \le t_1$.
 \cr}
\end{displaymath}
Then
\begin{displaymath}
V(t_1,x_1)\leq \int_0^{t_1} L(y_1(s), y_1^\prime(s))ds +
\varphi(y_1(t_1))
= \int_0^{h_1} L(x_1+s u_1, u_1)ds
\end{displaymath}
\begin{displaymath}
 \mbox{}+
\int_{s_{1}}^{t_{2}} L(y(s; t_2,x_2), y'(s; t_2,x_2))ds
+\varphi(y(t_{2}; t_2,x_2)).
\end{displaymath}
As $s_{1}=h_{1}-t_{1}+t_{2}\geq 0$ and $L\geq 0$, we obtain
\begin{displaymath}
V(t_1,x_1) \leq \int_0^{h_1} L(x_1+s u_1, u_1)ds +
V(t_2,x_2).
\end{displaymath}
Since $L$ is locally bounded, it follows from (\ref{dm301}) that
there exists a constant $M$, depending only on $L$, $t_0$, $x_0$,
$\delta$, and $r$, such that
\begin{displaymath}
V(t_1,x_1)-V(t_2,x_2) \leq Mh_1=M(|t_1-t_2| + |x_1-x_2|).
\end{displaymath}
Exchanging the roles of $(t_1,x_1)$ and $(t_2,x_2)$ we obtain that
$V$ in $M-$Lipschitz on $B ((t_0,x_0), \delta )$.
$\; \; \Box$

\vspace{ 3 mm}

When $L$ is discontinuous, different contingent inequalities have to
be introduced, which involve the function $L ^+(x,u)$ defined by
\begin{eqnarray}
L^+(x,u) &=& \limsup_{h \rightarrow 0+} \frac{1}{h} \inf \left\{
\int_{-h}^{0} L (y (s),y' (s))ds  : y (-h)=x-hu,\; y (0)=x\right\}
\nonumber \\
&=& \limsup_{h \rightarrow 0+} \frac{1}{h} \inf \left\{
\int_{0}^{h} L (y (s),-y' (s))ds  : y (0)=x ,\;   y (h)=x-hu\right\} .
\label{feb30}
\end{eqnarray}

\begin{Remark}\label{feb31}
{\rm
The function $L^+ (x,u)$ was introduced in  \cite{amar}. In that
paper it was shown that, if for some positive constants $D$, $d$
and $p>1$  we have
\begin{displaymath}
 \forall \; (x,u) \in \R^n \times \R^n ,
 \;\; d |u|^p  \leq L (x,u) \leq D (1 + |u|^p) ,
\end{displaymath}
then $L^+ (x, \cdot)$  is continuous for every $x\in \R^{n}$ and
convex for almost every $x\in\R^{n}$.

Clearly for all $u$ the function $L^{+} ( \cdot ,u)$ is
smaller than or equal to the upper semicontinuous envelope
of $L( \cdot ,u)$.
}
\end{Remark}

\begin{Proposition} \label{eq---}
If $L$ is locally bounded, then $ L^+(x,u)\geq  L(x,u)$
for all $ (x,u)\in \R^n\times\R^n$. Moreover, if $v_h\to u$ as $h\to
0+$, then
\begin{displaymath} L^+ (x,u) =\limsup_{h \rightarrow 0+}\frac{1}{h} \inf \left\{ \int_{-h}^{0} L (y ( s ),y' ( s )) d s
\;:\; y (-h)=x-hv_h  ,\;  y (0) = x\right\} .
\end{displaymath}
In particular $L^+=L$ when $L$
is continuous.
\end{Proposition}

{\bf Proof} --- \hspace{ 2 mm}
Let us fix $ (x,u)$ and $v_h$ as in the statement of the
proposition, and let $ \overline L (x,u)$ be the right hand side
of  the formula to be proved.
We want to show that $\overline L(x,u)\leq L^+(x,u)$.
For every $h>0$ let $\varepsilon_{h}= |v_h-u|$
and let $y_{h}$ be a minimizer of the problem
\begin{displaymath}
 \inf \left\{ \int_{-(1-\varepsilon_{h})h}^{0} L(y ( s ),y' ( s )) d s  :
y(-(1-\varepsilon_{h})h)=x-(1-\varepsilon_{h})h u  ,\;
 y (0)=x\right\} .
\end{displaymath}
For all $ s \in [-h,-(1-\varepsilon_{h})h]$ set
\begin{displaymath}
 y _h ( s) =x-(1-\varepsilon_{h})h u + (s+(1-\varepsilon_{h}) h)
 \frac{ v_h - (1-\varepsilon_{h}) u}{\varepsilon_{h}} .
\end{displaymath}
Then $y_h (-h)= x-h v_h$.
Since $ (v_h- (1-\varepsilon_{h}) u)/\varepsilon_{h}$ is bounded,
we deduce that for some $M>0$ and all $h>0$,
\begin{displaymath}
 \int_{-h}^{ 0} L (y_h (s),y_h' (s))ds \leq
\int_{-(1-\varepsilon_{h}) h}^{0}L (y_h
(s),y_h' (s))ds + \int_{-h}^{ -(1-\varepsilon_{h}) h} M ds .
\end{displaymath}
Dividing by $h$ and taking the upper limit
when $h \rightarrow 0+$
we get $\overline L (x,u) \leq L^+ (x,u)$. The opposite inequality can
be proved in the same way.

To prove that $L^+(x,u) \geq L(x,u)$,
for every $h>0$ let $y_{h}$ be a minimizer of the problem
\begin{displaymath}
\inf \left\{
\int_{-h}^{0} L (y ( s ),y' ( s )) d s  :  y (-h)=
x-hu ,\;  y (0)=x\right\} .
\end{displaymath}
Observe that $L^+(x,u)$ is
finite, because $L$ is locally bounded.
By \cite{ambr} there exist $k_h$ such that for some $p_h (s) \in \partial
_u L (y_h (s),y_h' (s))$
\begin{displaymath}
 \left\langle   p_h (s),y_h' (s)\right\rangle -L (y_h (s),y_h' (s))=k_h .
\end{displaymath}
We set $y_h (s)=y_h (-h)= x-h u$ for $s \leq -h$. Then for some $M
\geq 0$
\begin{displaymath}
 \int_{-1}^{0}L (y_h (s),y_h' (s))ds \leq  M .
\end{displaymath}
Arguing as in the proof of Proposition~\ref{minlip}, we deduce that
the sequence $k_h$ is bounded. By  \cite[Proof of Theorem
4.2]{ambr} the sequence $y_h'$ is bounded in $L^{ \infty}$.

Let $ \varepsilon >0$.
Since $L$ is lower semicontinuous and continuous in the second variable,
for all $h$ small enough
\begin{displaymath}
 L ( y_h (s), y_h' (s)) \geq L (x, y_h' (s)) - \varepsilon\;\; \mbox{ {\rm
 a.e.\ in}}\;\; [-h,0]
\end{displaymath}
(see Step 2 in the proof of Proposition \ref{cont}).
Integrating and using Jensen's inequality we get
\begin{displaymath}
\frac{1}{ h} \int_{-h}^{0 } L ( y_h (s), y_h' (s)) ds \geq
\frac{1}{h}\int_{-h}^{ 0}L (x, y_h' (s))ds - \varepsilon
\geq L (x,   u) - \varepsilon .
\end{displaymath}
Taking the upper limit when $h \rightarrow 0+ $ and $ \varepsilon
\rightarrow 0+$ we obtain
$L^+ (x,u) \geq L (x,u)$.

If $L$ is continuous, we have also $L^+\le L$ by Remark~\ref{feb31},
hence $L^+=L$.
$\; \; \Box$

\vspace{ 5 mm}

We generalize now a result obtained in  \cite{amar} under some additional
growth conditions.
\begin{Proposition} \label{ae}
Assume that $L$ is locally bounded.
Let $y\colon  [ 0,T] \mapsto  \R^n$ be a Lipschitz function. Then
$ 
 L (y (t),y' (t))= L^+ (y (t),y' (t))
$ 
for almost every $t \in [0,T]$.
\end{Proposition}

{\bf Proof} --- \hspace{ 2 mm}
We already know that $L^+\geq L$. To prove
the opposite inequality, we fix a function $y$ as in the statement of
the proposition.
Since $L$ is locally bounded, the function $t \mapsto
\psi (t):= \int_{0}^{t} L (y (s),y' (s))ds$ is absolutely continuous. Let
$t\in [0,T]$ be such that both $ \psi ' (t)$ and $y' (t)$ do exist and
$\psi ' (t) = L ( y (t),y' (t)) $. Since
$v_{h}=(y(t)-y(t-h))/h$ converges to $y'(t)$ as $h\to 0+$,
from Proposition~\ref{eq---} we obtain
\begin{displaymath}
L^+(y(t),y' (t)) \leq \lim_{h\to 0+} \frac{1}{h} \int_{t-h}^{t} L ( y
(s),y' (s))
ds = \psi ' (t) = L ( y (t),y'(t)) ,
\end{displaymath}
which concludes the proof.
$\; \; \Box$

\section{The Hamilton-Jacobi equation}

In this section we suppose that $L$ is continuous and we consider
three different notions of generalized solutions of
the Hamilton-Jacobi equation (\ref{H-J}). In Theorem~\ref{smaller}
we prove that
the value function $V$ defined by (\ref{value}) is the unique
lower semicontinuous solution of
the Hamilton-Jacobi equation (\ref{H-J}) in a suitable generalized
sense, expressed in terms of contingent inequalities.
Then we restrict our attention to
locally Lipschitz solutions, and we look for uniqueness results for
weaker (and more usual) notions of solutions. In Theorem~\ref{hju1}
we prove that
$V$ is the unique locally
Lipschitz function which satisfies the Hamilton-Jacobi equation
in the following generalized sense:
\begin{displaymath}
\forall \;(t,x) \in \R_+^{ \star }\times \R^n,\;\; \forall \;
(p_t,p_x)
\in \partial _- V (t,x),\;\; p_t+ H (x,-p_x)=0 .
\end{displaymath}
Finally, in Theorem~\ref{gg013}, we prove that $V$
is the unique locally Lipschitz viscosity solution of
(\ref{H-J}).

\begin{Theorem} \label{smaller}
Assume that $L$ is continuous.
Then $V$ is
lower semicontinuous on $\R_+ \times \R^n$
and solves the Hamilton-Jacobi equation (\ref{H-J}) in the following
sense:
\begin{equation} \label{cd22}
 \forall \;(t,x) \in \mbox{\rm dom}(V), \;t >0,\; \exists \; u
 \in\R^{n},\; \;
 D_{\uparrow }V (t,x) (-1, u) \leq -L (x,  u) ,
\end{equation}
 \begin{equation} \label{cd21}
 \forall \; (t,x) \in \mbox{\rm dom}(V), \;
 \forall \; u \in \R^n,\;\; D _{\downarrow }V (t,x) (1,-u) \leq L (x,u) .
\end{equation}

Let $W\colon  \R_+ \times \R^n \mapsto \R_+ \cup \{ + \infty \}$ be a
lower semicontinuous function which satisfies the initial condition
$W (0, \cdot )= \varphi$. If $W$ is a supersolution of the
Hamilton-Jacobi
equation (\ref{H-J}), in the sense that
\begin{equation} \label{febcd22}
 \forall \;(t,x) \in \mbox{\rm dom}(W), \;t >0,\; \exists \; u
 \in\R^{n},\; \;
 D_{\uparrow }W (t,x) (-1, u) \leq -L (x,  u) ,
\end{equation}
then $W\ge V$ on $\R_+ \times \R^n$. If $W$ is a subsolution
of the Hamilton-Jacobi
equation (\ref{H-J}), in the sense that
\begin{equation} \label{feb45}
 \forall \; (t,x) \in \mbox{\rm dom}(W), \;
 \forall \; u \in \R^n,\;\; D _{\downarrow }W (t,x) (1,-u) \leq L (x,u) ,
\end{equation}
then $W\le V$ on $\R_+ \times \R^n$.

Therefore $V$ is the unique non negative lower
semicontinuous function
which satisfies the initial condition $V(0, \cdot )= \varphi$
and solves the Hamilton-Jacobi
equation (\ref{H-J}) in the sense of (\ref{cd22}) and (\ref{cd21}).
\end{Theorem}

{\bf Proof} --- \hspace{ 2 mm}
The lower semicontinuity of $V$ is proved in Proposition~\ref{lsc}.
Condition (\ref{cd22}) follows from Proposition~\ref{cont}.
Let us prove (\ref{cd21}). Fix $(t,x) \in \mbox{\rm dom}(V)$ and
$u\in\R^n$.
Observe that for all $h>0$ and $v\in\R^n$ we have
\begin{displaymath}
 V (t+h, x-hv) - V(t,x) \leq \int_{0}^{h}L (x-sv,v)ds .
\end{displaymath}
Dividing by $h$ and taking the upper limit
when $h \rightarrow 0+$ and $v\to u$ we obtain (\ref{cd21}).

Let $W$ be as in the statement of the theorem. Assume
that $W$ is a supersolution, i.e., $W$ satisfies (\ref{febcd22}).
Let $ (t,x) \in \mbox{\rm dom}(W)$, $t>0$,
and $ \varepsilon >0$.
By (\ref{dirder}) and (\ref{febcd22}) there exist
$h_i \rightarrow 0+$ and $u_i \rightarrow u$ such that
\begin{displaymath}
 W (t-h_i,x +h_iu_i) \leq W (t,x) -h_i L (x,u) + \varepsilon h_i .
\end{displaymath}
Since $L$ is continuous, for all sufficiently large $i$ we have
\begin{displaymath}
 W (t-h_i,x +h_iu_i)+ \int_{0}^{h_i} L (x +
su_i,u_i)ds  \leq W (t,x) + 2 \varepsilon h_i .
\end{displaymath}
Consider the set $\cal A$ of all pairs $(\sigma, y)\in {] 0,t]} \times
W^{1,1}
(0,\sigma;\R^n)$ such that
\begin{displaymath}
y (0)=x\;\;\&\;\;  W (t-\sigma,y(\sigma))
+ \int_{0}^{\sigma} L(y (s),y' (s))ds  \leq W (t,x) + 2 \varepsilon
\sigma .
\end{displaymath}
The set $\cal A$ is nonempty by the first part of the proof.
We introduce the following partial order on $\cal A$: we write
$ (\sigma,y) << (\tau,z)$ if $\sigma\leq \tau$ and $y$ is the
restriction
of $z$ to $[0,\sigma]$.

We claim that for every totally ordered subset $\cal B$ of $\cal A$
there exists $
(\tau, z) \in {\cal A}$ such that
$(\sigma, y)  << (\tau,  z)$ for every
$(\sigma, y) \in {\cal B}$.
Indeed let
$$\tau = \sup_{ (\sigma,y) \in {\cal B}} \sigma$$
and consider a sequence $ (\sigma_i,y_i) \in {\cal B} $ with
$\sigma_i \rightarrow \tau$.
Define $z\colon {[0,\tau[} \mapsto \R^n$ by setting $z(s)=y_i (s)$
whenever
$s<\sigma_i$. As $\cal B$ is totally ordered, the function $z$ is well
defined
and, if $ (\sigma,y) \in {\cal B}$ with $\sigma<\tau$, then $y$
coincides with  the
restriction of $z$ to $[0,\sigma]$.

Since $W \geq 0$ we deduce that for some $c$ and for all $i$
\begin{displaymath}
 \int_{0}^{t_i} L (y_i(s),y_i' (s))ds \leq c .
\end{displaymath}
Set $y_i (s)=y_i (\sigma_i)$ for all $s \in {]\sigma_i,\tau]}$.
Since $L$ satisfies the coercivity assumption (\ref{T}), $y_i'$ are
equiintegrable on  $[0, \tau ]$. Consequently, $y_i$ are
equicontinuous
on  $[0, \tau ]$. So $z$ can be extended by continuity at  $ \tau $,
$z$ belongs
to $W^{1,1}(0,\tau;\R^{n})$ and $(\tau,z)$ belongs to $\cal A$
(recall that $W$ is lower semicontinuous). It is clear from the
construction that $(\sigma,y)<<(\tau,z)$ for every $ (\sigma,y) \in
{\cal
B}$.

By Zorn's lemma there exists a maximal
element $ (\sigma,y) \in {\cal A}$. We claim that $\sigma=t$.
Indeed, if not, then
acting as at the beginning of the proof, we construct $ (\tau,z)\in
{\cal A}$ with $\sigma<\tau\le t$ and
$(\sigma,y)<<(\tau,z)$, which contradicts the maximality. So
\begin{displaymath}
 V(t,x)\le W (0, y (t))+ \int_{0}^{t} L(y (s),y' (s))ds  \leq W (t,x)
+ 2
\varepsilon t .
\end{displaymath}
  The number $ \varepsilon >0$ being arbitrary, we conclude that
  $V(x,t)\le W(x,t)$.

 Suppose now that $W$ is a subsolution, i.e., $W$ satisfies
 (\ref{feb45}).
Let us fix $t>0$ and $x \in \R^n$, and
let $y$  be a minimizer of  (\ref{value}). Since $L$ is
continuous,  (\ref{bound}) holds true.
Thus $y' \in L^{ \infty } (0,t;\R^{n})$ by \cite{ambr}.
Consider a sequence of continuous functions $u_i\colon
[0,t] \mapsto \R^n$  which is bounded in $L^{ \infty } (0,t;\R^{n})$
and converges to $y'$ almost
everywhere in $[0,t]$,  and define
\begin{displaymath}
 y_i (s)=  y(t)  -\int_{s}^{t} u_i ( \tau )d \tau .
\end{displaymath}
Then $y_i$ converges to $y$ uniformly in $[0,t]$.
Fix  $i$ and set $ \psi (s)=W (s,y_i (t-s))$ for $0\le s \le t$.
By (\ref{feb45}) for every $s\in \mbox{\rm dom}( \psi)$, with $s <t$,
we have
\begin{equation}
 \label{coine}
D _{\uparrow } \psi  (s) (1) \leq
L (y_i (t-s),  u_i(t-s)).
\end{equation}

Consider the closed set
$$
K= {\cal E}pi ( \psi) \cup ({[t, +\infty [} \times \R) ,
$$
and the system
\begin{equation}\label{feb65}
\cases{
(\alpha' (s) ,  z' (s))=
(1, L (  y_i (t-s),  u_i (t-s))), &
\cr\cr
(\alpha  (0) , z (0))  = (0, \varphi (y (t))) . &
\cr}
\end{equation}
It has the unique solution
\begin{equation}\label{feb66}
( \alpha  (s),z (s)) =  \left( s, \varphi  (  y (t)) +
\int_{0}^{s} L (y_i (t- \tau ), u_i ( t-\tau )) d \tau
\right) .
\end{equation}
According to (\ref{star}), Theorem \ref{viability}, and
(\ref{coine}),
this solution is
viable in $K$, i.e.,  for all $s \in [ 0,t]$,
 $ ( \alpha  (s),z (s)) \in K$. Thus for all $s<t$
\begin{equation}\label{feb67}
W (s,y_i (t-s)) \leq  \varphi (y (t)) + \int_{0}^{t} L ( y_i (
\tau ),   u_i ( \tau ))d \tau.
\end{equation}
Since $L$ is continuous and $W$ is
lower semicontinuous, from the Lebesgue
Theorem we obtain
\begin{equation}\label{feb68}
 W (t,x) \leq \varphi (y  (t)) + \int_{0}^{t} L ( y ( \tau
),   y' ( \tau ))d \tau ,
\end{equation}
which gives $W (t,x) \leq V(t,x)$.
$\; \; \Box$

\begin{Remark}\label{feb22}
{\rm In the previous theorem, the comparison result for subsolutions
remains true, if we assume that (\ref{feb45}) holds only for
$t>0$, provided that $W$ satisfies also
condition (\ref{hj2}) of the next theorem. In this case (\ref{coine})
holds only for $s>0$, but we can modify the proof in the following way.
We fix $t_{i}\to 0+$ and $x_{i}\to y(t)$
such that $W(t_{i},x_{i})\to \varphi(y(t))$, and define
\begin{displaymath}
 y_i (s)=  x_{i}  -\int_{s}^{t-t_{i}} u_i ( \tau )d \tau ,
\end{displaymath}
for $0\le s\le t-t_{i}$
Consequently, $\psi(s)$ is defined
only for $t_{i}\le s\le t$,
and (\ref{coine}) holds only for $t_{i}\le s\le t$.
Then we replace $0$ by $t_{i}$ and
$\varphi(y(t))$ by $W(t_{i},x_{i})$ in (\ref{feb65})--(\ref{feb67}),
and we obtain (\ref{feb68}) as before.
}
\end{Remark}

\begin{Remark}\label{feb22a}
{\rm The following example shows that we can not remove the assumption  
$W \geq 0$ in Theorem \ref{smaller} even if $L$ does not depend on  $x$. Let
$n=1,\; L (x,u)= \frac{1}{2} |u|^2,\; \varphi (x)=0$. Then  $V (t,x)=0$ and
$H (x,p)= \frac{1}{2}|p|^2$. Let us consider the lower semicontinuous function
$W: \R_+ \times \R^n \mapsto \R$ defined by 
\begin{displaymath}
W (t,x) \;=\; \left\{ \begin{array}{ll}
0 & \mbox{ {\rm if }} tx<1 ,\\
- (tx-1)^{1/2} - x^3-xt^{-2} & \mbox{ {\rm if }} tx \geq 1 .
\end{array} \right. \end{displaymath}
By direct computation one checks that for  $tx \ne 1$ the function $W$ is 
differentiable and
\begin{displaymath}
 W_t (t,x) + H (x,-W_x (t,x)) =W _t (t,x) + \frac{1}{2}|W_x (t,x)|^2 \geq 0 ,
\end{displaymath}
which implies that
\begin{equation}
 \label{feb22b}
\exists \; u \in \R^n,\;\; D _{\uparrow }W (t,x) (-1,u) \leq -L (x,u).
\end{equation}
On the other hand for  $tx=1$ we have
\begin{displaymath}
\forall \; u > 1/t^2 ,\;\;
 D _{\uparrow }W (t,x) (-1,u)= - \infty ,\end{displaymath}
which implies  (\ref{feb22b}) also in this case. Therefore $W$ satisfies
(\ref{febcd22}), but  $W (t,x) <V (t,x)$ for  $tx \geq 1$.
}
\end{Remark}

The following proposition shows the equivalence between the notion of
subsolution considered in the previous theorem and a notion defined
by using subdifferentials.

\begin{Proposition} \label{feb60}
Assume that $L$ is continuous.
Let $W\colon  \R_+ \times \R^n \mapsto \R \cup \{ + \infty \}$ be a
lower semicontinuous function. Then the following conditions are
equivalent:
\begin{equation} \label{feb61}
 \forall \; (t,x) \in \mbox{\rm dom}(W), \;
 \forall \; u \in \R^n,\;\; D _{\downarrow }W (t,x) (1,-u) \leq L
 (x,u) ;
\end{equation}
\begin{equation}
 \label{feb61a}
 \forall \; (t,x) \in \mbox{\rm dom}(W), \;
 \forall \; u \in \R^n,\;\; D _{\uparrow }W (t,x) (1,-u) \leq L
 (x,u) ;
\end{equation} 
 \begin{equation} \label{feb62}
 \forall \; (t,x) \in \R_+ \times \R^n, \; \; \forall \; (p_t,p_x)
\in
\partial_- W (t,x), \;p_t + H (x,-p_x) \leq 0 .
\end{equation}
The equivalence remains true if $\mbox{\rm dom}(W)$ is replaced by
$\mbox{\rm dom}(W)\cap (\R_+^{ \star }\times \R^n)$ in (\ref{feb61}) and 
(\ref{feb61a})
and $\R_+\times \R^n$ is replaced by $\R_+^{ \star }\times
\R^n$ in (\ref{feb62}).
\end{Proposition}

{\bf Proof} --- \hspace{ 2 mm}
It is clear that (\ref{feb61}) implies (\ref{feb61a}).
Suppose that $W$ satisfies (\ref{feb61a}).
Then (\ref{feb62}) follows from (\ref{subd1}) and
 from the definition (\ref{H}) of $H$.

 Conversely, suppose that $W$ satisfies (\ref{feb62}).
 We claim that for all $(t,x) \in \mbox{\rm dom}(W)$
\begin{equation} \label{cd2}
 \forall \; u\in\R^{n},\;\; (1,-u, L (x,u) ) \in \overline {co}
\;T_{{\cal
E}pi (W)} (t,x, W (t,x)) .
\end{equation}
To prove this fact, let us fix $u \in \R^n$. Then
\begin{equation} \label{cd3}
 \forall \; (p_t,p_x) \in \partial_- W (t,x),\; \;p_t + \left\langle
-p_x,u\right\rangle -L (x,u)  \leq 0 .
\end{equation}
We want to prove that
\begin{equation} \label{cd4}
\forall\; (p_t,p_x,q) \in
\left[   T_{ {\cal E}pi (W)} (t,x,W(t,x))\right]^- ,\;\;
p_t + \left\langle   -p_x,u\right\rangle +qL (x,u)  \leq 0 .
\end{equation}
When $q<0$ this inequality follows from (\ref{subd}) and (\ref{cd3}).
By Lemma  \ref{Rockafellar}, if $ (0,0,0) \ne (p_t,p_x,0) \in
\left[ T_{{\cal E}pi (W)} (t,x, W (t,x)) \right] ^-$,
then for some $ (t_i,x_i) \rightarrow (t,x)$
and
$(p^i_{t},p^i_{x},q_i) \in \left[   T_{{\cal E}pi (W)} (t_i,x_i, W
(t_i,x_i)) \right]^-$,
with $q_i <0$, we have $(p^i_{t},p^i_{x},q_i) \rightarrow
(p_t,p_x,0)$.
So
\begin{displaymath}
 p^i_{t} + \left\langle   -p^i_{x},u\right\rangle +q_iL (x_i,u)  \leq
0 .
\end{displaymath}
Taking the limit we get
$p_t + \left\langle   -p_x,u\right\rangle  \leq 0$, which concludes
the proof of (\ref{cd4}).

By the separation theorem, (\ref{cd2}) follows from
(\ref{cd4}). Thus for all $ (t,x,r)
\in {\cal E}pi (W) $
\begin{displaymath}
 (1,-u,L (x,u)) \in \overline {co} \;T_{{\cal
E}pi (W)} (t,x, r) .
\end{displaymath}
{}From Theorem \ref{corn} and continuity of $L$ we deduce that for
all $
(t,x) \in \mbox{\rm dom}(W),$
\begin{displaymath}
 (1,-u,L (x,u)) \in C_{ {\cal E}pi(W)}  (t,x,W (t,x)) .
\end{displaymath}
Fix $\varepsilon >0$.
Then it is not difficult to check that
\begin{displaymath}
\forall \; u \in \R^n,\;\;  (1,-u,L (x,u) +\varepsilon)
\in \mbox{\rm Int} \left(
C_{ {\cal E}pi(W)}
 (t,x,W (t,x)) \right) .
\end{displaymath}
By \cite[Proposition 13, p. 425]{AE} this yields
\begin{equation} \label{B+}
 \limsup_{h \rightarrow 0+, \,v \rightarrow u} \frac{W (t +h, x-hv)-W
(t,x)}{h}  \leq L (x,u) +\varepsilon .
\end{equation}
As $\varepsilon>0$ is arbitrary, we obtain (\ref{feb61}).
$\; \; \Box$

\vspace{ 5 mm}

We recall that for a locally Lipschitz function
$ \varphi \colon  \R^n \mapsto \R$
the contingent derivatives coincides with the Dini
derivatives:
\begin{displaymath}
\forall \; u \in \R^n,\;\; D _{\uparrow } \varphi  (x) (u) =
d_- \varphi(x;u):=
\liminf_{h \rightarrow 0+} \frac{ \varphi (x+hu)- \varphi
(x)}{h} < +\infty .
\end{displaymath}
\begin{displaymath}
\forall \; u \in \R^n,\;\; D _{\downarrow } \varphi  (x) (u) =
d_+\varphi(x;u):=
\limsup_{h \rightarrow 0+} \frac{ \varphi (x+hu)- \varphi
(x)}{h} < +\infty .
\end{displaymath}
By (\ref{star}) this implies that
\begin{equation}
\label{ddd5}
(p,0)\in\left[T_{{\cal E}pi(\varphi)}(x,\varphi(x))\right]^-\;\;
\Longrightarrow \;\;
p=0.
\end{equation}

\vspace{ 5 mm}

\begin{Theorem} \label{hju1}
Assume that $L$ is continuous.
Then $V$
is locally Lipschitz on $\R_+^{ \star } \times \R^n$
and solves the Hamilton-Jacobi equation
(\ref{H-J}) in the following sense:
\begin{equation} \label{hj1}
\forall \;(t,x) \in \R_+^{ \star }\times \R^n,\;\; \forall \;
(p_t,p_x)
\in \partial _- V (t,x),\;\; p_t+ H (x,-p_x)=0.
\end{equation}

Let $W \colon \R_+^{ \star } \times \R^n \mapsto \R_+$ be a locally
Lipschitz function which satisfies the initial condition
\begin{equation} \label{hj2}
\forall\; x\in\R^{n},\;\;
 \liminf_{h \rightarrow 0+, \;y \rightarrow x} W (h,y) = \varphi(x) .
\end{equation}
If $W$ is a supersolution of  the Hamilton-Jacobi equation
(\ref{H-J}), in the sense that
\begin{equation} \label{feb01}
\forall \;(t,x) \in \R_+^{ \star }\times \R^n,\;\; \forall \;
(p_t,p_x)
\in \partial _- W (t,x),\;\; p_t+ H (x,-p_x)\ge 0,
\end{equation}
then $W\ge V$ on $\R_+^{ \star }\times \R^n$.
If $W$ is a subsolution of  the Hamilton-Jacobi equation
(\ref{H-J}), in the sense that
\begin{equation} \label{feb02}
\forall \;(t,x) \in \R_+^{ \star }\times \R^n,\;\; \forall \;
(p_t,p_x)
\in \partial _- W (t,x),\;\; p_t+ H (x,-p_x)\le 0,
\end{equation}
then $W\le V$ on $\R_+^{ \star }\times \R^n$.

Therefore the restriction of $V$ to $\R_+^{ \star } \times \R^n$ is
the
unique non negative locally Lipschitz function which satisfies the initial
condition
(\ref{hj2}) and solves the Hamilton-Jacobi
equation (\ref{H-J}) in the sense of (\ref{hj1}).
\end{Theorem}

{\bf Proof} --- \hspace{ 2 mm}
The fact that $V$ is locally Lipschitz on $\R_+^{ \star } \times \R^n$
is proved in Corollary~\ref{vlip}. Property (\ref{hj1})
follows from Theorem~\ref{smaller}. Indeed, (\ref{cd22}) and
(\ref{subd1}), together with the definition (\ref{H}) of $H$, imply
that
\begin{displaymath}
\forall \;(t,x) \in \R_+\times \R^n,\;\; \forall \;
(p_t,p_x)
\in \partial _- V (t,x),\;\; p_t+ H (x,-p_x)\ge 0;
\end{displaymath}
the opposite inequality follows from (\ref{cd21}) and
Proposition~\ref{feb60}.
The initial condition (\ref{hj2}) for $V$
is proved in Proposition~\ref{lsc}.

Let $W$ be as in the statement of the theorem. Assume that
$W$ is a
viscosity supersolution, i.e., $W$ satisfies (\ref{feb01}).
Define
\begin{displaymath}
  G (x) = \left\{  (-1,u,-L (x,u) - \rho ) :  \rho  \geq 0,\; u \in
\R^n \right\} ,
\end{displaymath}
and fix $t>0$ and  $x \in \R^n$. Since $W$ is Lipschitz around
$(t,x)$,
 $ \partial _- W$ is locally bounded. By
the coercivity assumption (\ref{T}) there exist $R>0$ and $\delta >0$
such
that
\begin{eqnarray*}
 &\forall \; (s,z) \in B ((t,x), \delta ),\;\; \forall \; (p_t,p_x)
\in
\partial _-W (s,z), \;\;\forall \; |u| \geq R, &\\
&p_t + \left\langle
-p_x,u\right\rangle -L (z,u) <0 .&
\end{eqnarray*}
This, (\ref{feb01}),
 and the separation theorem imply that for every
 $(s,z) \in B ((t,x), \delta )$
\begin{displaymath}
G (z) \cap \left(
\left\{ -1 \right\} \times B (0,R) \times [-m,0] \right) \cap
\left( \overline{co}\, T_{ {\cal E}pi (W)} (s,z,W (s,z)) \right) \ne
\emptyset ,
\end{displaymath}
where $m = \max\{L (z,u):  (z,u) \in {B (x, \delta ) \times B (0,R)}
\}$.
The above holds true also with $W (s,z)$ replaced by any $ r  \geq W
(s,z)$.
{}From Theorem \ref{viability} we obtain
\begin{displaymath}
 G (x) \cap \left(  \{-1\}
\times B (0,R) \times [-m,0] \right) \cap
\left(  T_{ {\cal E}pi (W)} (t,x,W (t,x)) \right) \ne \emptyset ,
\end{displaymath}
which is equivalent to
\begin{displaymath}
 \exists \; u \in B (0,R),\;  D _{\uparrow }W (t,x) (-1,u) \leq -L
(x,u) .
\end{displaymath}
{}From Theorem \ref{smaller} we deduce that $W \geq V$ on
$\R_+^{\star}
\times \R^n$.

If $W$ is a subsolution, i.e., $W$ satisfies (\ref{feb02}), then
$W\le V$ on $\R_+^{\star} \times \R^n$
by Theorem~\ref{smaller} and Remark~\ref{feb22}.
$\; \; \Box$

\begin{Theorem}\label{gg013}
Assume that $L$ is continuous.
Then $V$
is locally Lipschitz on $\R_+^{ \star } \times \R^n$
and solves the Hamilton-Jacobi equation
(\ref{H-J}) in the viscosity sense, i.e.,
\begin{equation} \label{ddd10}
 \forall \; (t,x) \in \R_+^{ \star }\times \R^n, \;\; \forall \;
(p_t,p_x)
\in \partial _- V (t,x),\;\; p_t+ H (x,-p_x) \geq 0 ,
\end{equation}
\begin{equation} \label{ddd11}
  \forall \; (t,x) \in \R_+^{ \star }\times \R^n, \;\; \forall \;
(p_t,p_x)
\in \partial _+ V (t,x),\;\; p_t+ H (x,-p_x) \leq 0 .
\end{equation}

Let $W \colon \R_+^{ \star } \times \R^n \mapsto \R_+$ be a locally
Lipschitz function which satisfies the initial condition
(\ref{hj2}).
If $W$ is a viscosity supersolution of the Hamilton-Jacobi
equation (\ref{H-J}), i.e.,
\begin{equation} \label{febddd10}
 \forall \; (t,x) \in \R_+^{ \star }\times \R^n, \;\; \forall \;
(p_t,p_x)
\in \partial _- W (t,x),\;\; p_t+ H (x,-p_x) \geq 0 ,
\end{equation}
then $W\ge V$ on $\R_+^{ \star }\times \R^n$. If
$W$ is a viscosity subsolution of the Hamilton-Jacobi
equation (\ref{H-J}), i.e.,
\begin{equation} \label{febddd11}
  \forall \; (t,x) \in \R_+^{ \star }\times \R^n, \;\; \forall \;
(p_t,p_x)
\in \partial _+ W (t,x),\;\; p_t+ H (x,-p_x) \leq 0 ,
\end{equation}
then $W\le V$ on $\R_+^{ \star }\times \R^n$.

Therefore the restriction of  $V$ to $\R_+^{ \star } \times \R^n$ is
the
unique non negative locally Lipschitz viscosity solution of the Hamilton-Jacobi
equation (\ref{H-J}) which satisfies the initial condition
(\ref{hj2}).
\end{Theorem}

{\bf Proof} --- \hspace{ 2 mm}
Let us prove (\ref{ddd11}).
Fix $(t,x) \in \mbox{\rm dom}(V)$, $t>0$, and
$u\in\R^n$.
Observe that for all small $h>0$ we have
\begin{displaymath}
 V(t,x) \leq V (t-h, x+hu)+ \int_{0}^{h}L (x+su,u)ds ,
\end{displaymath}
hence
\begin{displaymath}
V (t-h, x+hu)) - V(t,x) \geq -\int_{0}^{h}L (x+su,u)ds .
\end{displaymath}
Dividing by $h$ and taking the upper limit
when $h \rightarrow 0+$ we obtain
\begin{displaymath}
D_{\downarrow } V (t,x) (-1,u) \geq  -L(x,u) ,
\end{displaymath}
which implies (\ref{ddd11}) by (\ref{supd1}) and by the definition
(\ref{H}) of $H$.
For the other properties of
$V$ see Theorem~\ref{hju1}.

Let $W$ be as in the statement of the theorem. If $W$ is a
viscosity supersolution, i.e., $W$ satisfies (\ref{febddd10}), then
$W\ge V$ on $\R_+^{ \star }\times \R^n$ by Theorem~\ref{hju1}.

Assume now that $W$ is a viscosity subsolution, i.e., $W$ satisfies
(\ref{febddd11}).
Properties (\ref{subd}) and (\ref{ddd5}) imply that, for all $t>0$ and
$x\in\R^{n}$,
\begin{displaymath}
 \left[ T_{{\cal E}pi(-W)}(t,x,-W(t,x))\right]^-=
 \bigcup_{ \lambda \geq 0}^{} \lambda \left(
-\partial _+ W (t,x),-1 \right) .
\end{displaymath}
Using the separation theorem, from (\ref{febddd11}) we obtain
\begin{displaymath}
 \forall \; u \in \R^n,\;\; (-1,u, L (x,u)) \in \overline{ \mbox{{\rm
co}}}\; T_{{\cal E}pi(-W)}(t,x,-W(t,x)) .
\end{displaymath}
The above holds true also with $-W (t,x)$ replaced by any $r \geq -W (t,x)$.
So, by Theorem \ref{corn}, for all $t>0$ and
$x\in\R^{n}$,
\begin{displaymath}
\forall \; u \in \R^n,\;\; (-1,u, L (x,u)) \in
\,T_{{\cal E}pi(-W)}(t,x,-W(t,x)),
\end{displaymath}
and therefore by (\ref{star}) and (\ref{subd1})
\begin{displaymath}
\forall \; u \in \R^n,\;\; D _{\downarrow }W (t,x)
(-1,u)
\geq -L (x,u) .
\end{displaymath}
Fix  $t>0, x \in \R^n$ and let $y$ be a solution to (\ref{value}). By
\cite{ambr} it is Lipschitz. Consider  $ (t_i,x_i) \rightarrow  (0+,y
(t))$
such that $\displaystyle
\lim_{i \rightarrow \infty }W (t_i,x_i)= \varphi (y (t))$. Set
\begin{displaymath}
 y_i (s)= y (s)+x_i -y (t-t_i) .
\end{displaymath}
The function $ \psi  (s):= W (t-s,y_i (s))$
is locally Lipschitz on $[0,t[$. Fix $0 \leq s<t$ such that $ \psi '
(s)$
and  $y' (s)$ do exist.
Using the fact that $W$ is locally Lipschitz we get
\begin{displaymath}
 \psi ' (s)= D _{\downarrow }W (t-s,y _i (s)) (-1,y' (s)) \geq -L (y_i
(s),y'_i
(s)) .
\end{displaymath}
Consequently for every $0 \leq s < t$
\begin{displaymath}
 W (t-s, y_i (s)) - W (t,y_i (0)) = \psi (s) - \psi  (0) \geq
-\int_{0}^{s}
L (y_i (\tau ),y'_i
(\tau ))d \tau ,
\end{displaymath}
and thus
\begin{displaymath}
 W (t_i, x_i) + \int_{0}^{t} L (y_i ( \tau ),y' ( \tau )) d\tau
\geq W (t,y_i (0)) .
\end{displaymath}
Passing to the  limit when  $i \rightarrow \infty $ and using
continuity
of
$L$ we deduce that  $W (t,x) \leq V (t,x)$.
$\; \; \Box$

\section{The case of a discontinuous Lagrangian}

In this section we do not assume that $L$ is continuous. We can still
prove (Theorem~\ref{last}) that the value function $V$ defined by
(\ref{value}) is the unique non negative
lower semicontinuous solution of
the Hamilton-Jacobi equation (\ref{H-J}), but now we have to consider
a weaker notion of generalized solution, which involves a contingent
inequality for the function
$L^+$ introduced in (\ref{feb30}). To prove uniqueness in the smaller
class of locally Lipschitz functions we can use
an even weaker notion of
solution, where the contingent inequality (\ref{cd25}) for
supersolutions is
replaced by a viscosity inequality (Theorem~\ref{ddd050}).

\begin{Theorem} \label{last}
Assume that $L$ is locally bounded.
Then $V$ is
lower semicontinuous on $\R_+ \times \R^n$
and solves the Hamilton-Jacobi equation (\ref{H-J}) in the following
sense:
\begin{equation} \label{cd25}
 \forall \; (t,x) \in \mbox{\rm dom}(V),\; t>0, \;
 \exists\; u\in\R^{n}, \;\;
 D _{\uparrow }V (t,x) (-1,u) \leq -L (x,u) ,
\end{equation}
\begin{equation} \label{cd26}
\forall \; (t,x) \in \mbox{\rm dom}(V),\;
 \forall \; u\in\R^{n},\;\; D_{\downarrow } V (t,x) (1,-u) \leq  L^+
(x,u) .
\end{equation}

Let $W\colon  \R_+ \times \R^n \mapsto \R_+ \cup \{ + \infty \}$ be a
lower semicontinuous function which satisfies the initial condition
$W (0, \cdot )= \varphi$. If $W$ is a supersolution of the
Hamilton-Jacobi
equation (\ref{H-J}), in the sense that
\begin{equation} \label{feb04}
 \forall \; (t,x) \in \mbox{\rm dom}(W),\; t>0, \;
 \exists\; u\in\R^{n}, \;\;
 D _{\uparrow }W (t,x) (-1,u) \leq -L (x,u) ,
\end{equation}
then $W\ge V$ on $\R_+ \times \R^n$. If $W$ is a subsolution
of the Hamilton-Jacobi
equation (\ref{H-J}), in the sense that
\begin{equation} \label{feb05}
\forall \; (t,x) \in \mbox{\rm dom}(W),\;
 \forall \; u\in\R^{n},\;\; D _{\downarrow } W (t,x) (1,-u) \leq  L^+
(x,u) ,
\end{equation}
then $W\le V$ on $\R_+ \times \R^n$.

Therefore $V$ is the
unique non negative lower semicontinuous  function
which satisfies the initial condition $V(0, \cdot )= \varphi$
and solves the Hamilton-Jacobi
equation (\ref{H-J}) in the sense of (\ref{cd25}) and (\ref{cd26}).
\end{Theorem}

To prove the theorem, we need the following approximation lemma.

\begin{Lemma}\label{approx}
There exists a sequence of
continuous functions $L_k\colon \R^{n} \times \R^n \mapsto  \R_+ $,
converging pointwise to $L$, such that, for all $k$,
$L_k (x, \cdot )$ is convex,
$L_k \leq L_{k+1} \leq L$, and
\begin{equation}\label{feb20}
\forall \; x\in \R^{n},
\;\forall \; u\in \R^{n},
\;\; L_k (x,u) \geq  \Theta  ( u ) .
\end{equation}
For every $k$ let $V_k$ be the value function of problem
(\ref{value}) with $L$ replaced by $L_k$. Then
$V_k$ converge to  $V$ pointwise.
\end{Lemma}

{\bf Proof} --- \hspace{ 2 mm} The proof of the existence of a
sequence
$L_{k}$ with the required properties can be found in
\cite[Lemma 2.2.3]{butt}. It is clear that the sequence $V_k$
is nondecreasing, so it is pointwise convergent, and that $V_{k}\le
V$.
We want to prove that $V\le \lim_{k} V_{k}$. Let us fix $(t,x)$ in
$\R_+\times \R^n$ such that
\begin{equation} \label{ff11}
\lim_{k\to\infty} V_{k}(x,t) < +\infty .
\end{equation}
Let $y_k$ be a solution to the problem
\begin{displaymath}
\inf \left\{ \int_{0}^{t} L_k (y (s),y' (s))ds \;+\; \varphi (y
(t)) : \; y(0)=x \right\} .
\end{displaymath}
We deduce from (\ref{ff11}) and from the coercivity assumption
(\ref{T}) that $y'_k$ are equiintegrable. Hence $y_k$ are
equicontinuous.
Taking a
subsequence and keeping the same notations we may assume that $y_k$
converges
uniformly to some $y$ and $y_k'$ converges weakly in $L^1 (0,t;\R^n)$
to
$y'$. Fix $i$. Since for every $k \geq i$ we have
\begin{displaymath}
 \varphi  (y_k (t)) + \int_{0}^{t} L _i (y_k (s),y_k ' (s))ds \leq
\varphi  (y_k (t)) + \int_{0}^{t} L _k (y_k (s),y_k ' (s))ds = V_{k}
(t,x) ,
\end{displaymath}
taking the limit when $ k \rightarrow \infty $ we get
\begin{displaymath}
 \varphi  (y (t)) + \int_{0}^{t} L _i (y(s),y ' (s))ds \leq
\lim_{k \rightarrow \infty } V_k (t,x) \leq  V (t,x) .
\end{displaymath}
Taking the limit when $i \rightarrow \infty $ and using Fatou's lemma
we
deduce that
\begin{displaymath}
V (t,x) \leq  \varphi  (y (t)) + \int_{0}^{t} L (y(s),y ' (s))ds \leq
\lim_{k \rightarrow \infty } V_k (t,x) \leq  V
(t,x) .
\end{displaymath}
Thus
$V_k$ converges pointwise to  $V$.
$\; \; \Box$

\vspace{ 3 mm}

{\bf Proof of Theorem \ref{last}} --- \hspace{ 2 mm}
The lower semicontinuity of $V$ is proved in Proposition \ref{lsc}.
Condition (\ref{cd25}) follows from Proposition
\ref{cont}. Let us prove (\ref{cd26}).
Fix $(t,x) \in \mbox{\rm dom}(V)$ and $u\in\R^{n}$.
By Proposition~\ref{eq---} we have
\begin{displaymath}
L^{+} (x,u)= \limsup_{ \textstyle {h \rightarrow 0+ \atop v \rightarrow u}}
\frac{1}{h} \inf \left\{
\int_{-h}^{0} L (y ( s ), y' ( s )) d s  :  y(-h)=x-h v, \; y (0)=x\right\} .
\end{displaymath}
Since for every absolutely continuous function $y$ satisfying
$y(-h)=x-h v$ and $y(0)=x$ we have
\begin{displaymath}
 V (t+h,x-h v) \leq V(t,x) + \int_{-h}^{0} L (y (s),y' (s)) ds ,
\end{displaymath}
we deduce that
\begin{displaymath}
 V (t+h,x-h v ) - V (t,x)
\end{displaymath}
\begin{displaymath}
 \leq  \inf \left\{  \int_{-h}^{0} L (y(s),y' (s)) ds  :
  y (0)=x,\; y (-h)=x-h v \right\} .
\end{displaymath}
Dividing by $h$ and taking the upper limit as $h\to 0+$ and $v\to u$,
we obtain (\ref{cd26}).

Let $W$ be as in the statement of the theorem. Assume
that $W$ is a supersolution, in the sense that (\ref{feb04}) is
satisfied. Let $L_{k}$ and $V_{k}$ be the continuous Lagrangians and
the corresponding value functions given by Lemma~\ref{approx}.
As $L_{k}\le L$, the function $W$ is a supersolution for the problem
relative to
the continuous Lagrangian $L_{k}$. Therefore $W\ge V_{k}$
on $\R_+ \times \R^n$ by Theorem~\ref{smaller}.
Since $V_{k}$ converges to $V$ pointwise, we
conclude that $W\le V$ on $\R_+ \times \R^n$.

Assume now that $W$ is a subsolution, i.e.,
$W$ satisfies (\ref {feb05}). Fix $ t>0,\; x\in \R^n$ and
let $y$
be a minimizer of  (\ref{value}). It is Lipschitz continuous by
\cite{ambr}.
Set $ \psi(s)=W (s,y (t-s))$. Thus for almost all $s \in [ 0,t]$,
\begin{equation} \label{coine1}
D _{\uparrow } \psi  (s) (1) \leq D _{\downarrow } W (s,y(t-s))
(1,-y'(t-s)) \leq  L^{+}(y(t-s),  y' (t-s)) .
\end{equation}
By Proposition \ref{ae} we have
\begin{equation}\label{feb70}
L^{+}(y(t-s ), y' (t- s ))= L(y(t-s ), y' (t- s ))
\;\; \mbox {\ a.e.\ in\ } [0,t].
\end{equation}

Since $y$ is Lipschitz and  $L$ is locally bounded, the same argument
used for (\ref{coine1})
implies that there exists a constant  $M $ such that
\begin{equation}\label{ff20}
\forall \; s \in {[0,t[}, \;\;D _{\uparrow }
\psi (s) (1) \leq M.
\end{equation}
Define the closed set-valued map $s \leadsto P (s)$ by
$$ P (s)=  W (s, y (t-s)) + \R_+,\; \forall \;  s \in {[0,t[} \; \;\;
\&\;\;\; P(s) = \R,\;\; \forall \;  s \geq t .$$
Using (\ref{ff20}), (\ref{star}), and
Theorem \ref{viability}, we deduce that
for every $s_{0}\in [0,t]$ and for every $z_{0}\in P(s_{0})$
there exists an $M$--Lipschitz function
$z\colon{[s_{0},t]}\mapsto \R^{n}$ such that $z(s_{0})=z_{0}$
and $z (s) \in P (s)$ for every $s\in [s_{0},t]$. This yields that
$P$ is left
absolutely continuous on $[0,t]$, i.e.,
for any $\varepsilon >0$ and for any compact
subset $K \subset  \R^n$, there exists $\delta > 0$
such  that for any subdivision
$0\leq t_1 < \tau _ 1 \leq  \ldots  t_m < \tau _m \leq t$ with
$\sum_{i}(\tau _i -t_i)  \leq \delta$ we have
$\sum_{i} h(P( t_i) \cap K, P( \tau _i)) \leq \varepsilon$,
where $h$ is the Hausdorff semidistance:
$h(A,B):= \displaystyle\sup_{a\in\,A} d(a,B)$.

Consider the viability problem
\begin{equation}\label{feb80}
\cases{
 z' (s)  = L (  y (t-s),  y' (t-s)),&
 \cr
z (0)  = \varphi (y (t)), &
\cr
z (s) \in P (s) .&
\cr}
\end{equation}
According to the measurable viability theorem
\cite[Theorem 4.2]{FPR92} it has a (viable)
solution by (\ref{coine1}) and (\ref{feb70}).
But this solution is given by
$$
z (s) =   \varphi  (y (t)) + \int_{0}^{s} L(y(t-
\tau ), y' (t- \tau )) d \tau .$$
Thus
$$W (t,x) \leq  \varphi (y (t)) + \int_{0}^{t} L ( y (
\tau ),   y' ( \tau ))d \tau =V (t,x) , $$
and $W \leq V$ on $\R_+ \times \R^n$.
$\; \; \Box$

\vspace{ 5 mm}

The proof shows that the comparison result for supersolutions
in Theorem~\ref{last} remains
true even if we drop the assumption that $L$ is locally bounded.
Therefore Proposition~\ref{cont} and Theorem~\ref{last} imply the
following corollary.

\begin{Corollary}\label{gg1}
Assume that $L$ satisfies (\ref {bound}). Then $V$
is the smallest non negative lower
semicontinuous function satisfying
the initial condition $V (0, \cdot )= \varphi$ and the
contingent inequality (\ref{feb04}).
\end{Corollary}

\begin{Remark}\label{feb75}
{\rm
In Theorem~\ref{last} the comparison result for subsolutions remains
true, if we assume only that
(\ref{feb05}) holds for $t>0$, provided that $W$ satisfies condition
(\ref{ddd150}) of the next theorem.
In this case (\ref{ff20}) holds only for $s>0$,
but, to obtain the result, it is enough to replace $0$ by $h>0$
and $\varphi(t)$ by $W(h, y(t-h))$ in (\ref{feb80}).
}
\end{Remark}

We consider now the uniqueness in the class of locally Lipschitz
solutions.

\begin{Theorem}\label{ddd050}
Assume that $L$ is locally bounded.
Then $V$
is locally Lipschitz on $\R_+^{ \star } \times \R^n$
and solves the Hamilton-Jacobi equation
(\ref{H-J}) in the following sense:
\begin{equation} \label{hj3}
\forall \;(t,x) \in \R_+^{ \star }\times \R^n,\;\; \forall \;
(p_t,p_x)
\in \partial _- V (t,x),\;\; p_t+ H (x,-p_x) \geq 0,
\end{equation}
\begin{equation}
 \label{hj4}
\forall \;(t,x) \in \R_+^{ \star }\times \R^n,\;\; \forall \;u \in
\R^n,
\;\; D_{\downarrow }V (t,x) (1,-u) \leq L^+ (x,u) .
\end{equation}

Let $W \colon \R_+^{ \star } \times \R^n \mapsto \R_+ $ be a locally
Lipschitz function which satisfies the initial condition (\ref{hj2})
together with
\begin{equation} \label{ddd150}
\forall \; x\in \R^{n}, \; \; \forall \; \lambda>0, \;\;
\lim_{\scriptstyle h \rightarrow 0+, \;y \rightarrow x
\atop \scriptstyle |y-x| \le \lambda h} W (h,y) = \varphi  (x).
\end{equation}
If $W$ is a supersolution of  the Hamilton-Jacobi equation
(\ref{H-J}), in the sense that
\begin{equation} \label{feb07}
\forall \;(t,x) \in \R_+^{ \star }\times \R^n,\;\; \forall \;
(p_t,p_x)
\in \partial _- W (t,x),\;\; p_t+ H (x,-p_x) \geq 0,
\end{equation}
then $W\ge V$ on $\R_+^{ \star }\times \R^n$.
If $W$ is a subsolution of  the Hamilton-Jacobi equation
(\ref{H-J}), in the sense that
\begin{equation} \label{feb08}
\forall \;(t,x) \in \R_+^{ \star }\times \R^n,\;\; \forall \;u \in
\R^n,
\;\; D_{\downarrow }W (t,x) (1,-u) \leq L^+(x,u) ,
\end{equation}
then $W\le V$ on $\R_+^{ \star }\times \R^n$.

Therefore the restriction of $V$ to $\R_+^{ \star } \times \R^n$ is
the
unique non negative locally Lipschitz function which satisfies the initial
conditions
(\ref{hj2}) and (\ref{ddd150}) and solves the Hamilton-Jacobi
equation (\ref{H-J}) in the sense of (\ref{hj3}) and (\ref{hj4}).
\end{Theorem}

{\bf Proof} --- \hspace{ 2 mm}
The fact that $V$ is locally Lipschitz on $\R_+^{ \star } \times \R^n$
is proved in Corollary \ref{vlip}. Conditions (\ref{hj2}) and
(\ref{ddd150})  for $V$
follow from Proposition \ref{lsc}. Condition (\ref{hj3}) is proved in
Proposition \ref{cont}, while (\ref{hj4}) follows from Theorem
\ref{last}.

Let $W$ be as in the statement of the theorem. If
$W$ is a supersolution, i.e., $W$ satisfies (\ref{feb07}),
then we can prove that $W \geq V$, arguing as in the
proof of Theorem~\ref{hju1} (with Theorem~\ref{smaller} replaced by
Theorem~\ref{last}).
If $W$ is a subsolution, i.e., $W$
satisfies
(\ref{feb08}), then $W\le V$ by Theorem~\ref{last} and
Remark~\ref{feb75}.
$\; \; \Box$

\end{document}